\documentclass[]{article}
\usepackage{amssymb}

\begin{document}
\title{ Solutions of Analytical Systems of Partial Differential Equations}

\date{}

\author{ Kostadin Tren\v{c}evski\\
Institute of Mathematics, Faculty of Natural Sci. and Math.,\\
Ss. Cyril and Methodius University in Skopje, Macedonia\\
e-mail: kostadin.trencevski@gmail.com,  \quad kostatre@pmf.ukim.mk} 

\maketitle
\begin{abstract} In this paper are examined 
general classes of linear and non-linear analytical  systems 
of partial differential equations. Indeed the  integrability 
conditions are found and if they are satisfied, the 
solutions are given as functional series in a neighborhood of 
a given point $(x=0)$.  
\end{abstract}



\section{Introduction} 
This paper is a continuation of the papers \cite{2,3,4,5}, and 
we will give a brief view of them. 
\par 
In the paper \cite{2}  it  was  found  a  formula  for  the $k$-th
covariant derivative. Further that formula  was  generalized  for
$k\in \mathbb{R}$. Specially, if $k=-1$ it 
yields to a  general  solution  for  a 
system of linear differential equations \cite{3}. 
In the paper \cite{4} is given an application of \cite{3} for solving 
the Frenet equations. In \cite{4} two main theorems are proved. 
The first theorem gives the solution of analytical 
non-homogeneous linear system of differential equations 
of order $k$ of $n$ equations and $n$ unknown functions. 
The second theorem gives the solution of non-linear 
analytical system of differential equations (of the first 
order) of $n$ equations and $n$ unknown functions. 
\par 
In this paper we will prove two main theorems, 
considering linear and non-linear 
systems of partial differential equations. 
Without loss of generality, we will find the required solutions 
in a neighborhood of the point $(0, \cdots ,0).$ The 
author is not known to similar results by any other author. 
\par
The results of this paper have 
applications in the differential geometry \cite{5}, in 
studying the non-linear  connections  \cite{1}.  For  example  the 
integrability conditions  in  this  paper  are  just  
vanishing of the curvature tensor of the corresponding 
connections.  If   the   systems   of   partial 
differential equations considering in this paper are  tensor 
equations, then the obtained solutions also have tensor 
character. 

\section{Homogeneous system of linear partial differential equations} 
Let us consider the following system
\begin{equation}
\partial y_{r}/\partial x_{u} + \sum^{n}_{s=1}f_{rsu}y_{s} = 0\qquad 
(1\le r\le n, 1\le u\le k)
\label{2.1}
\end{equation}
\noindent of unknown functions $y_{1},\cdots  
,y_{n}$ of $k$ variables $x_{1},\cdots 
,x_{k}$ and $f_{rsu}$ are given analytical functions 
of $x_{1},\cdots ,x_{k}$, regular in a neighborhood of 
$(0,\cdots ,0)$. In order 
to  consider  the  integrability  conditions,  we  introduce  the
following functions
\begin{equation}
R_{tsuv} = \partial f_{tsv}/\partial x_{u} - \partial 
f_{tsu}/\partial x_{v} + \sum^{n}_{p=1} f_{tpu}f_{psv} - 
\sum^{n}_{p=1} f_{tpv}f_{psu}.
\label{2.2}
\end{equation}
$$
(1\le u,v\le k, 1\le t,s\le n)
$$
\noindent If (\ref{2.1}) is an integrable system for 
arbitrary initial conditions, then using that
$$
{\partial \over \partial x_{v}} {\partial y_{r}\over \partial x_{u}} = 
{\partial \over \partial x_{u}} {\partial y_{r}\over \partial x_{v}}
$$
\noindent and the system (\ref{2.1}), it is easy to obtain that
\begin{equation}
R_{tsuv} \equiv  0
\label{2.3}
\end{equation}
\noindent for $1\le u,v\le k$ and $1\le t,s\le $n. Conversely, it is known  
that  if  (\ref{2.3})
are satisfied, then the system (\ref{2.1}) is integrable. Indeed  this
statement also follows from the theorem 2.1.
\par
{\bf Theorem 2.1.}  {\it Let  the  system}  (\ref{2.1})  {\it 
with the initial conditions} $y_{s}(0,\cdots ,0)=C_{s}$ 
   $(1\le s\le n)$ 
{\it is given  and   the
integrability conditions} (\ref{2.3}) {\it are satisfied.  Then  there  exist
functions} 
$P^{<w_{1},...,w_{k}>}_{ts}(x_{1},\cdots 
,x_{k})$, 
 $w_{1},\cdots ,w_{k}\in \mathbb{ N}_{0}   
{\it and}   1\le t,s\le n,$ 
{\it such that}
\begin{equation}
P^{<0,...,0>}_{ts} = \delta _{ts};
\label{2.4a}
\end{equation}
$$
P^{<w_{1},...,w_{u}+1,...,w_{k}>}_{ts} = {\partial 
\over \partial x_{u}} 
P^{<w_{1},...,w_{k}>}_{ts} +
$$
\begin{equation}
+\sum^{n}_{p=1}f_{tpu}P^{<w_{1},...,w_{k}>}_{ps} 
\label{2.4b}
\end{equation}
\noindent  {\it  and  the  solution  of}  (\ref{2.1})  {\it  in  a 
neighborhood of} $(0, \cdots ,0)$ {\it is given by}
$$
y_{r}=\sum^{n}_{s=1} \sum^{\infty }_{w_{1}=0} 
\sum^{\infty }_{w_{2}=0} \cdots 
\sum^{\infty }_{w_{k}=0} {(-x_{1})^{w_{1}}\over w_{1}!} 
{(-x_{2})^{w_{2}}\over 
w_{2}!} \cdots 
 {(-x_{k})^{w_{k}}\over w_{k}!}\cdot 
$$
\begin{equation}
\cdot P^{<w_{1},...,w_{k}>}_{rs}C_{s}. \qquad (1\le r\le n)
\label{2.5}
\end{equation}
\noindent {\it This solution is unique with the given 
initial conditions in a neighborhood of} $(0,\cdots ,0)$. 
\par
Proof. Let us suppose that the system (\ref{2.1}) is given and the
integrability conditions (\ref{2.3}) are satisfied. In order  to  prove
that there exist functions 
$P^{<w_{1},...,w_{k}>}_{ts}(x_{1},\cdots 
,x_{k}) (w_{1},\cdots 
,w_{k}\in \mathbb{ N}_{0},
1\le t,s\le n)$
\par
\noindent such  that (\ref{2.4a}  and (\ref{2.4b})  are  satisfied,  it  is
sufficient to prove that
\begin{equation}
P^{<w_{1},...,w^{(2)}_{u}+1,...,w^{(1)}_{v}+1,...,w_{k}>}_{ts} = 
P^{<w_{1},...,w^{(1)}_{u}+1,...,w^{(2)}_{v}+1,...,w_{k}>}_{ts}
\label{2.6}
\end{equation}
\noindent for each $t,s\in \{1,\cdots ,n\}$ and $u,v\in \{1,\cdots 
,k\}, u\neq v$,  where  the 
notations $"(1)"$ and $"(2)"$ show the order of the two increased 
indices. Indeed 
$$
P^{<w_{1},...,w^{(2)}_{u}+1,...,w^{(1)}_{v}+1,...,w_{k}>}_{ts} 
= {\partial \over \partial 
x_{u}} P^{<w_{1},...,w_{v}+1,...,w_{k}>}_{ts} +
$$
$$
+ \sum^{n}_{p=1} f_{tpr}P^{<w_{1},...,w_{v}+1,...,w_{k}>}_{ps} =
$$
$$
= {\partial \over \partial x_{u}} \Bigl[ 
{\partial \over 
\partial x_{v}} P^{<w_{1},...,w_{k}>}_{ts} + \sum^{n}_{q=1} 
f_{tqv}P^{<w_{1},...,w_{k}>}_{qs} \Bigr] +
$$
$$
+ \sum^{n}_{p=1} f_{tpu} \Bigl[ {\partial \over \partial 
x_{v}} P^{<w_{1},...,w_{k}>}_{ps} + \sum^{n}_{a=1} 
f_{pav}P^{<w_{1},...,w_{k}>}_{as} \Bigr]
$$
\noindent and similarly
$$
P^{<w_{1},...,w^{(1)}_{u}+1,...,w^{(2)}_{v}+1,...,w_{k}>}_{ts} =
$$
$$
= {\partial \over \partial x_{v}} \Bigl[ 
{\partial \over 
\partial x_{u}} P^{<w_{1},...,w_{k}>}_{ts} + \sum^{n}_{q=1} 
f_{tqu}P^{<w_{1},...,w_{k}>}_{qs} \Bigr] +
$$
$$
+ \sum^{n}_{p=1} f_{tpv} \Bigl[ {\partial \over \partial 
x_{u}} P^{<w_{1},...,w_{k}>}_{ps} + \sum^{n}_{a=1} 
f_{pau}P^{<w_{1},...,w_{k}>}_{as} \Bigr] .
$$
\noindent Hence we obtain 
$$
P^{<w_{1},...,w^{(2)}_{u}+1,...,w^{(1)}_{v}+1,...,w_{k}>}_{ts} - 
P^{<w_{1},...,w^{(1)}_{u}+1,...,w^{(2)}_{v}+1,...,w_{k}>}_{ts} =
$$
$$
= \sum^{n}_{q=1} R_{tquv}P^{<w_{1},...,w_{k}>}_{qs},
$$
\noindent and (\ref{2.6}) is satisfied because $R_{tquv}\equiv 0.$
\par
Now we should prove that the functions $(y_{r})$ from (\ref{2.5})  
satisfy the system (\ref{2.1}). 
First we prove the uniform convergence of the right side of 
(\ref{2.5}) in a neighborhood of $(0,\cdots ,0)$. We can 
consider analytical functions of complex variables. Suppose 
that $x=(x_{1},\cdots ,x_{k})$ is sufficiently close to 
$(0,\cdots ,0)$ such that all functions $\{ f_{rsu}\}$ 
are regular in the disc $D_{x}=\{ z=(z_{1},\cdots ,z_{k}): 
\mid z-x \mid <\rho \}$ and $0\in D_{x}$. Hence $\mid x/\rho \mid <1.$ 
Obviously, all functions $P^{<w_{1},\cdots ,w_{k>}}_{ts}$ are 
regular in $D_{x}$. In order to find an 
estimation of $P^{<w_{1},\cdots ,w_{k}>}_{ts}$ from (\ref{2.4a})  
and (\ref{2.4b}), some additional results should be given. Let 
$D_{x_{u}}  (1\le u\le k)$ be an operator defined by 
\begin{equation} 
D_{x_{u}}(y_{r})=\partial y_{r}/\partial x_{u} + \sum^{n}_{s=1}f_{rsu}y_{s} \qquad 
(1\le r\le n). 
\label{2.7}
\end{equation}
\noindent If the integrability conditions (\ref{2.3}) are 
satisfied, then similarly to the result in \cite{2}, it holds 
$$
(D^{w_{1}}_{x_{1}} \circ D^{w_{2}}_{x_{2}}\circ  \cdots  \circ 
D^{w_{k}}_{x_{k}})(y_{r}) = \sum^{w_{1}}_{m_{1}=0} \cdots 
 \sum^{w_{k}}_{m_{k}=0} \sum^{n}_{s=1} P^{<m_{1},...,m_{k}>}_{rs}\cdot 
$$
\begin{equation}
\cdot {\partial ^{w_{1}-m_{1}+...+w_{k}-m_{k}}y_{s}\over \partial 
x^{w_{1}-m_{1}}_{1}\partial x^{w_{2}-m_{2}}_{2}...\partial 
x^{w_{k}-m_{k}}_{k}}\cdot {w_{1}! \cdots 
w_{k}!\over m_{1}!\cdots 
m_{k}!(w_{1}-m_{1})!\cdots 
(w_{k}-m_{k})!}. 
\label{2.8}
\end{equation}
\noindent Since $P^{<w_{1},...,w_{k}>}_{rj}=(D^{w_{1}}_{x_{1}} \circ  \cdots  
\circ D^{w_{k}}_{x_{k}})\delta _{rj}$ for fixed $j$, by putting 
$y_{r}=P^{<1,...,1>}_{rj}$, we obtain
\par
$$
P^{<w_{1}+1,...,w_{k}+1>}_{rj} = \sum^{w_{1}}_{m_{1}=0} \cdots 
 \sum^{w_{k}}_{m_{k}=0} \sum^{n}_{s=1} P^{<m_{1},...,m_{k}>}_{rs}\cdot 
$$
$$
\cdot \Bigl[ {\partial ^{w_{1}-m_{1}+...+w_{k}-m_{k}}\over \partial 
x^{w_{1}-m_{1}}_{1}...\partial x^{w_{k}-m_{k}}_{k}} P^{<1,...,1>}_{sj} 
\Bigr] {w_{1}! \cdots 
w_{k}!\over m_{1}!\cdots 
m_{k}!(w_{1}-m_{1})!\cdots 
(w_{k}-m_{k})!}.
$$
\noindent This equality is suitable for estimation of $\mid 
P^{<w_{1},...,w_{k}>}_{rj}\mid $.
\par
If $Q^{<w_{1},...,w_{k}>}_{rj}=P^{<w_{1},...,w_{k}>}_{rj}/(w_{1}!\cdots 
w_{k}!)$, then
\par
$$
Q^{<w_{1}+1,...,w_{k}+1>}_{rj} = {1\over (w_{1}+1)\cdots (w_{k}+1)} 
\sum^{w_{1}}_{m_{1}=0} \cdots 
 \sum^{w_{k}}_{m_{k}=0} \sum^{n}_{s=1} 
$$
\begin{equation}
Q^{<m_{1},...,m_{k}>}_{rs}
\cdot {\partial ^{w_{1}-m_{1}+...+w_{k}-m_{k}}\over \partial 
x^{w_{1}-m_{1}}_{1}...\partial 
x^{w_{k}-m_{k}}_{k}} P^{<1,...,1>}_{sj} {1\over 
(w_{1}-m_{1})!\cdots 
(w_{k}-m_{k})!}.
\label{2.9}
\end{equation}
\noindent According to the Cauchy integral formula, it holds
\par
\begin{equation}
\max _{s,j} \mid {\partial ^{w_{1}-m_{1}+...+w_{k}-m_{k}}\over \partial 
x^{w_{1}-m_{1}}_{1}...\partial x^{w_{k}-m_{k}}_{k}} P^{<1,...,1>}_{sj}\mid  
\le  {M\cdot (w_{1}-m_{1})! \cdots (w_{k}-m_{k})! 
\over \rho ^{w_{1}-m_{1}+...+w_{k}-m_{k}}},
\label{2.10}
\end{equation}
\noindent where $M$ depends (continuously) only on $x_{1},\cdots 
,x_{k}$. Let 
$$A^{<w_{1},...,w_{k}>}_{r} = \max _{j} \mid Q^{<w_{1},...,w_{k}>}_{rj}\mid .$$ 
\noindent Then (\ref{2.9}) and (\ref{2.10}) imply
\par
$$
A^{<w_{1}+1,...,w_{k}+1>}_{r} \le {1\over (w_{1}+1)\cdots (w_{k}+1)} 
\sum^{w_{1}}_{m_{1}=0} \cdots 
 \sum^{w_{k}}_{m_{k}=0} A^{<m_{1},...,m_{k}>}_{r}\cdot $$
$$
\cdot {nM\over \rho ^{w_{1}-m_{1}+...+w_{k}-m_{k}}}.
$$
\noindent Now  if  $\rho$  is  sufficient  small  such  that 
$nM\rho ^{k}\le 1,$ then
\par
$$A^{<w_{1}+1,...,w_{k}+1>}_{r} \rho ^{(w_{1}+1)+...+(w_{k}+1)} \le $$
\par
$$
\le  {1\over (w_{1}+1)\cdots (w_{k}+1)} \sum^{w_{1}}_{m_{1}=0} \cdots 
 \sum^{w_{k}}_{m_{k}=0} A^{<m_{1},...,m_{k}>}_{r}\rho ^{m_{1}+...+m_{k}}.
$$
\noindent Moreover we can suppose that instead of (\ref{2.10}) it holds
\par
$$
\max _{s,j} \mid {\partial ^{w_{1}-m_{1}+...+w_{k}-m_{k}}\over \partial 
x^{w_{1}-m_{1}}_{1}...\partial x^{w_{k}-m_{k}}_{k}} P^{<a_{1},\cdots 
,a_{k}>}_{sj}\mid  \le  {M(w_{1}-m_{1})!\cdots (w_{k}-m_{k})! 
\over \rho ^{w_{1}-m_{1}+...+w_{k}-m_{k}}}
$$
\noindent for each $a_{1},\cdots 
,a_{k}\in \{0,1\}$, and $\rho$ is such that $nM\rho ^{u}\le 1$  
for $1\le u\le $k.
Now by induction of $k$ it is easy to verify that
\par
$$
A^{<m_{1},...,m_{k}>}_{r}\rho ^{m_{1}+...+m_{k}} \le  1.
$$
\noindent Thus
\par
$$
\mid {1\over w_{1}!\cdots w_{k}!} P^{<w_{1},...,w_{k}>}_{rj}(x_{1},\cdots 
,x_{k})\mid  \le  {1\over \rho ^{m_{1}+...+m_{k}}},
$$
\noindent and we have uniform convergence in  (\ref{2.5})  for $\mid 
x_{1}\left|\matrix{,\cdots 
,}\right|x_{k}\mid \le (1-\epsilon ) \rho $.
According to the  Weierstrass  theorem,  the  functions $y_{r}$  are
regular in a neighborhood of $(0,\cdots 
,0)$  and  it  is  allowed  to
differentiate them by parts.
\par
If $C_{1},\cdots 
,C_{n}$  are  arbitrary  constants,
using (\ref{2.5}) and (\ref{2.4b}) we obtain
$$
\partial y_{r}/\partial x_{u} =
\sum^{n}_{s=1} 
\sum^{}_{w_{1},...,w_{k}\in \mathbb{ N_{0}}} 
{(-x_{1})^{w_{1}}\over w_{1}!}\cdots 
(-1){(-x_{u})^{w_{u}-1}\over (w_{u}-1)!}\cdots 
{(-x_{k})^{w_{k}}\over w_{k}!}\cdot 
$$
$$
\cdot P^{<w_{1},...,w_{k}>}_{rs}C_{s} +$$ 
$$+\sum^{n}_{s=1} 
\sum^{}_{w_{1},...,w_{k}\in \mathbb{ N_{0}}} 
{(-x_{1})^{w_{1}}\over w_{1}!}\cdots 
{(-x_{k})^{w_{k}}\over w_{k}!}\cdot {\partial \over \partial x_{u}} 
P^{<w_{1},...,w_{k}>}_{rs}C_{s}
$$
$$
= -\sum^{n}_{s=1} 
\sum^{}_{w_{1},...,w_{k}\in \mathbb{ N_{0}}} 
{(-x_{1})^{w_{1}}\over w_{1}!}\cdots 
{(-x_{u})^{w_{u}}\over w_{u}!}\cdots 
{(-x_{k})^{w_{k}}\over w_{k}!}\cdot 
$$
$$
\cdot P^{<w_{1},...,w_{u+1},...,w_{k}>}_{rs}C_{s} +
$$
$$
+\sum^{n}_{s=1} 
\sum^{}_{w_{1},...,w_{k}\in \mathbb{ N_{0}}} 
{(-x_{1})^{w_{1}}\over w_{1}!}\cdots 
{(-x_{k})^{w_{k}}\over w_{k}!}\cdot {\partial \over \partial x_{u}} 
P^{<w_{1},...,w_{k}>}_{rs}C_{s} =
$$
$$
= -\sum^{n}_{s=1} \sum^{}_{w_{1},...,w_{k}\in \mathbb{ N_{0}}} 
{(-x_{1})^{w_{1}}\over w_{1}!}\cdots 
{(-x_{k})^{w_{k}}\over w_{k}!}\cdot 
\Bigl[ P^{<w_{1},...,w_{u}+1,...,w_{k}>}_{rs} -
$$
$$
- {\partial \over \partial x_{u}} P^{<w_{1},...,w_{k}>}
_{rs} \Bigr] C_{s} =
$$
$$
= -\sum^{n}_{s=1} \sum^{}_{w_{1},...,w_{k}\in \mathbb{ N_{0}}} 
{(-x_{1})^{w_{1}}\over w_{1}!}\cdots 
{(-x_{k})^{w_{k}}\over w_{k}!} \sum^{n}_{p=1} 
f_{rpu}P^{<w_{1},...,w_{k}>}_{ps}C_{s} =
$$
$$
= - \sum^{n}_{p=1} f_{rpu}\Bigl[ \sum^{n}_{s=1} 
\sum^{}_{w_{1},...,w_{k}\in \mathbb{ N_{0}}} 
{(-x_{1})^{w_{1}}\over w_{1}!}\cdots 
{(-x_{k})^{w_{k}}\over w_{k}!} P^{<w_{1},...,w_{k}>}_{ps} C_
{s} \Bigr] 
$$
$$
= - \sum^{n}_{p=1} f_{rpu}y_{p},
$$
\noindent i.e. (\ref{2.1}) is satisfied. Moreover, using (\ref{2.4a}) we obtain
\par
$$
y_{r}(0,\cdots 
,0) = P^{<0,...,0>}_{rs} C_{s} = 
\delta _{rs}C_{s} = C_{r}. 
$$
\par 
To the end of the proof we have only to prove the uniqueness 
of the solution (\ref{2.1}).  Since  the  system  (\ref{2.1}) is 
linear, it is sufficient to prove that $y_{s}(0,\cdots ,0)=
C_{s}=0$  $(1\le s\le n)$  implies $y_{s}=0$   
$(1\le s\le n)$. Since the functions $f_{rsu}$ are 
analytical, each solution of (\ref{2.1}) is analytical. Using 
that $y_{s}(0,\cdots ,0)=0$, it follows  from  (\ref{2.1})  that 
the first partial derivatives of $y_{s}$ vanish at $(0,
\cdots ,0)$. By successive partial differentiations of 
(\ref{2.1}), all partial derivatives of $y_{s}$ vanish 
at $(0,\cdots ,0)$. Hence $y_{s}(x_{1},\cdots ,x_{k})=0$ in 
a neighborhood of $(0,\cdots ,0)$. \qquad $\Box$ 

{\bf Remark 1.} If one (or more) equations in the system (\ref{2.1}) is (are) omitted, then 
we may consider the whole system with $nk$ equations, where 
the corresponding functions $f_{rsu}$ can be considered as parameters. These parametric functions 
appear in the integrability conditions and also in the solution, and it leads to the parametric solution of the initial system.   

\section{Non-linear system of partial differential equations}
Let us consider the following non-linear system  of  partial
differential equations
$$
\partial y_{r}/\partial x_{u} + F_{ru}(x_{1},\cdots 
,x_{k},y_{1},\cdots 
,y_{n}) = 0 \qquad (1\le r\le n, 1\le u\le k).
$$
\noindent We suppose  that $F_{ru}(x_{1},\cdots 
,x_{k},y_{1},\cdots 
,y_{n})$  can  be  written  in  a
Laurent's series, i.e.
\par
\begin{equation}
\partial y_{r}/\partial x_{u}+ \sum^{}_{i_{1},...,i_{n}\in \mathbb{ Z}} 
f_{ri_{1}...i_{n}u}(x_{1},\cdots 
,x_{k})y^{i_{1}}_{1} y^{i_{2}}_{2}\cdots y^{i_{n}}_{n} = 0
\label{3.1}
\end{equation}
$$
1\le r\le n, 1\le u\le k
$$
\noindent where $f_{ri_{1}...i_{n}u}$ are analytical 
functions. Moreover, suppose that there exists a 
neighborhood U of $(0,\cdots ,0)$ such that all functions 
$f_{ri_{1}...i_{n}u}$ are regular in U. Let W is such that 
the Laurent's series in (\ref{3.1}) converge for 
$(y_{1},\cdots ,y_{n})\in W$ and $(x_{1},\cdots ,x_{k})\in U$. 
Before we consider  the  integrability  conditions  and  the
solution of the system (\ref{3.1}), we will introduce some notations.
\par
If $f_{ri_{1}...i_{n}u} (1\le r\le n,1\le u\le k, i_{1},\cdots 
,i_{n}\in \mathbb{ Z})$ are given functions
of $x_{1},\cdots 
,x_{k}$, then we define new functions $h_{i_{1}...i_{n}j_{1}...j_{n}u}$  and
$$R_{i_{1}...i_{n}j_{1}...j_{n}uv} (i_{1},\cdots 
,i_{n},j_{1},\cdots 
,j_{n}\in \mathbb{ Z},\quad  1\le u,v\le k).$$ 
\noindent First define 
\begin{equation}
h_{i_{1}...i_{n}j_{1}...j_{n}u} = \sum^{n}_{s=1} 
i_{s}f_{s(j_{1}-i_{1})...(j_{s}-i_{s}+1)...(j_{n}-i_{n})u}. 
\label{3.2}
\end{equation}
\noindent Now we will prove the convergence of the series 
$$
\sum^{}_{t_{1},...,t_{n}\in \mathbb{ Z}} 
h_{t_{1}...t_{n}j_{1}...j_{n}v}h_{i_{1}...i_{n}t_{1}...t_{n}u}. 
$$
\noindent According to the definition (\ref{3.2}), it is sufficient to 
prove the convergence of the series 
$$
\sum^{}_{t_{1},...,t_{n}\in \mathbb{ Z}} 
f_{p(j_{1}-t_{1})...(j_{p}-t_{p}+1)...(j_{n}-t_{n})v}
\cdot f_{s(t_{1}-i_{1})...(t_{s}-i_{s}+1)...(t_{n}-i_{n})u}. 
$$
\noindent Indeed, it converges 
because that is the coefficient in front of 
$$
z_{1}^{j_{1}-i_{1}}\cdots z_{s}^{j_{s}-i_{s}+1}\cdots 
z_{p}^{j_{p}-i_{p}+1}\cdots z_{n}^{j_{n}-i_{n}} 
$$
\noindent of the product of the Laurent's series
$$
\sum^{}_{t_{1},...,t_{n}\in \mathbb{ Z}} 
f_{p(j_{1}-t_{1})...(j_{p}-t_{p}+1)...(j_{n}-t_{n})v}\cdot 
z_{1}^{j_{1}-t_{1}}\cdots z_{p}^{j_{p}-t_{p}+1}\cdots
z_{n}^{j_{n}-t_{n}}\hbox{and} 
$$ 
$$
\sum^{}_{t_{1},...,t_{n}\in \mathbb{ Z}} 
f_{s(t_{1}-i_{1})...(t_{s}-i_{s}+1)...(t_{n}-i_{n})u}\cdot 
z_{1}^{t_{1}-i_{1}}\cdots z_{s}^{t_{s}-i_{s}+1}\cdots
z_{n}^{t_{n}-i_{n}},
$$
which are convergent for $(z_{1},\cdots ,z_{n})\in W$. 
Now we can define 
$$
R_{i_{1}...i_{n}j_{1}...j_{n}uv} = {\partial \over \partial x_{u}} 
h_{i_{1}...i_{n}j_{1}...j_{n}v} - {\partial \over \partial x_{v}} 
h_{i_{1}...i_{n}j_{1}...j_{n}u} +
$$
$$
+ \sum^{}_{t_{1},...,t_{n}\in \mathbb{ Z}} 
h_{t_{1}...t_{n}j_{1}...j_{n}v}h_{i_{1}...i_{n}t_{1}...t_{n}u} -
$$
\begin{equation}
- \sum^{}_{t_{1},...,t_{n}\in \mathbb{ Z}} 
h_{t_{1}...t_{n}j_{1}...j_{n}u}h_{i_{1}...i_{n}t_{1}...t_{n}v}.
\label{3.3}
\end{equation}
\noindent Note that the series 
$$
\sum^{}_{t_{1},...,t_{n}\in \mathbb{ Z}} 
h_{i_{1}...i_{n}j_{1}...j_{n}u}y_{1}^{j_{1}}\cdots 
y_{n}^{j_{n}}\qquad \hbox{and} \qquad 
R_{i_{1}...i_{n}j_{1}...j_{n}uv}y_{1}^{j_{1}}\cdots 
y_{n}^{j_{n}} 
$$
\noindent converge for $(y_{1}, \cdots ,y_{n})\in W$ and 
$(x_{1}, \cdots ,x_{k})\in U$. 
In order to simplify the notations, sometimes we will  denote  by
the Greek indices $\alpha ,\beta ,\gamma ,\cdots $
a set of $n$ integer indices $i_{1}\cdots 
i_{n}; j_{1}\cdots 
j_{n}; \cdots $ We will denote by $\{r\}$ the set of
$n$ indices $0\cdots 
010\cdots 
0$   where  the  unit  appears  at  the $r$-th
place. Now $\alpha +\beta$ and $\alpha -\beta$ are defined by 
$$
i_{1}\cdots 
i_{n} \pm  j_{1}\cdots 
j_{n} = (i_{1}\pm j_{1})(i_{2}\pm j_{2})\cdots 
(i_{n}\pm j_{n}).
$$
\par
{\bf Theorem 3.1.} {\it The quantities} $h_{\alpha \beta u}$ 
{\it and} 
$R_{\alpha \beta uv}$  {\it satisfy  the
following properties}:
\begin{equation}
h_{(\alpha +\beta )\gamma u} = h_{\alpha (\gamma -\beta )u} + h_{\beta 
(\gamma -\alpha )u},
\label{3.4}
\end{equation}
\begin{equation}
R_{(\alpha +\beta )\gamma uv} = R_{\alpha (\gamma -\beta )uv} + R_{\beta 
(\gamma -\alpha )uv},
\label{3.5}
\end{equation}
\begin{equation}
h_{\alpha \beta u} 
= \sum^{n}_{s=1} 
i_{s}h_{\{s\}(\beta -\alpha +\{s\})u},
\label{3.6}
\end{equation}
\begin{equation}
R_{\alpha \beta uv} 
= \sum^{n}_{s=1} 
i_{s}R_{\{s\}(\beta -\alpha +\{s\})uv},
\label{3.7}
\end{equation}
\noindent {\it where} $\alpha =i_{1}\cdots i_{n}. $
\par
Proof. Using the definition (\ref{3.2}) we obtain
$$
h_{\alpha (\gamma -\beta )u} + h_{\beta (\gamma -\alpha )u} =
$$
$$
= h_{i_{1}...i_{n}(t_{1}-j_{1})...(t_{n}-j_{n})u} + 
h_{j_{1}...j_{n}(t_{1}-i_{1})...(t_{n}-i_{n})u} =
$$
$$
= \sum^{n}_{s=1} i_{s}f_{s(t_{1}-j_{1}-i_{1})...(t_{n}-j_{n}-i_{n})u} + 
\sum^{n}_{s=1} j_{s}f_{s(t_{1}-i_{1}-j_{1})...(t_{n}-i_{n}-j_{n})u}=
$$
$$
= \sum^{n}_{s=1} 
(i_{s}+j_{s})f_{s(t_{1}-(i_{1}+j_{1}))...(t_{n}-(i_{n}+j_{n}))u} =
$$
$$
= h_{(i_{1}+j_{1})...(i_{n}+j_{n})t_{1}...t_{n}u} 
= h_{(\alpha +\beta )\gamma 
u},
$$
\par
\noindent and the identity (\ref{3.4}) is proved. 
\par
From the definition of $R_{\lambda \mu uv}$, i.e.
$$
R_{\lambda \mu uv} = {\partial \over \partial x_{u}} h_{\lambda \mu v} - 
{\partial \over \partial x_{v}} h_{\lambda \mu u} + \sum^{}_{\delta } 
h_{\delta \mu v}h_{\lambda \delta u} - \sum^{}_{\delta } h_{\delta \mu 
u}h_{\lambda \delta v}
$$
\noindent and the identity (\ref{3.4}) we obtain
$$
R_{(\alpha +\beta )\gamma uv} =$$
$$= {\partial \over \partial x_{u}}h_{\alpha 
(\gamma -\beta )v} + {\partial \over \partial x_{u}}
h_{\beta (\gamma -\alpha 
)v}- {\partial \over \partial x_{v}}
h_{\alpha (\gamma -\beta )u}- {\partial 
\over \partial x_{v}}h_{\beta (\gamma -\alpha )u} 
$$
$$
+ \sum^{}_{\delta } h_{\delta \gamma v}
(h_{\alpha (\delta -\beta )u}+h_{\beta 
(\delta -\alpha )u}) - \sum^{}_{\delta }
 h_{\delta \gamma u}(h_{\alpha (\delta 
-\beta )v}+h_{\beta (\delta -\alpha )v}) =
$$
$$
= {\partial \over \partial x_{u}} h_{\alpha 
(\gamma -\beta )v} - {\partial \over \partial x_{v}} h_{\alpha 
(\gamma -\beta )u} 
+ {\partial \over 
\partial x_{u}}h_{\beta (\gamma -\alpha )v} - {\partial \over \partial 
x_{v}}h_{\beta (\gamma -\alpha )u} + 
$$
$$
+ \sum^{}_{\delta } (h_{(\delta -\beta )
(\gamma -\beta )v} + h_{\beta (\gamma 
-\delta +\beta )v})h_{\alpha (\delta -\beta )u}
$$
$$
+ \sum^{}_{\delta } (h_{(\delta -\alpha )
(\gamma -\alpha )v} + h_{\alpha 
(\gamma -\delta +\alpha )v})h_{\beta (\delta -\alpha )u}
$$
$$
- \sum^{}_{\delta } (h_{(\delta -\beta )(\gamma 
-\beta )u} + h_{\beta (\gamma 
-\delta +\beta )u})h_{\alpha (\delta -\beta )v}
$$
$$
- \sum^{}_{\delta } (h_{(\delta -\alpha )(\gamma -\alpha )u} + h_{\alpha 
(\gamma -\delta +\alpha )u})h_{\beta (\delta -\alpha )v} =
$$
$$
= R_{\alpha (\gamma -\beta )uv} + R_{\beta (\gamma -\alpha )uv} + 
\sum^{}_{\delta } h_{\beta (\gamma -\delta 
+\beta )v}h_{\alpha (\delta -\beta 
)u} +
$$
$$
+ \sum^{}_{\delta } h_{\alpha (\gamma -\delta 
+\alpha )v}h_{\beta (\delta 
-\alpha )u} - 
\sum^{}_{\delta } h_{\beta (\gamma -\delta +\beta )u}h_{\alpha 
(\delta -\beta )v} -$$
$$- \sum^{}_{\delta } h_{\alpha (\gamma -\delta +\alpha 
)u}h_{\beta (\delta -\alpha )v} 
= R_{\alpha (\gamma -\beta )uv} + R_{\beta (\gamma -\alpha )uv}
$$
\par
\noindent because 
$$\sum^{}_{\delta } h_{\beta (\gamma -\delta +\beta 
)v}h_{\alpha (\delta -\beta )u} = \sum^{}_{\delta } 
h_{\alpha (\gamma -\delta 
+\alpha )u}h_{\beta (\delta -\alpha )v}$$
\par
\noindent and 
$$\sum^{}_{\delta } h_{\alpha (\gamma -\delta +\alpha 
)v}h_{\beta (\delta -\alpha )u} = \sum^{}_{\delta } 
h_{\beta (\gamma -\delta 
+\beta )u}h_{\alpha (\delta -\beta )v}.$$
\par
\noindent Hence the identity (\ref{3.5}) is proved.
\par
Finally, (\ref{3.6}) and (\ref{3.7}) are direct  consequences  of  (\ref{3.4})
and (\ref{3.5}). Indeed, using (\ref{3.4}) and (\ref{3.5}) one can verify  that  if
(\ref{3.6}) and (\ref{3.7}) hold for the set of indices $i_{1}\cdots 
i_{n}$,  then  they
also hold for  the  set  of  indices $i_{1}\cdots 
(i_{s}\pm 1)\cdots 
i_{n}$  for  each
$s\in \{1,\cdots 
,n\}$.
\par
We notice that (\ref{3.6}) can be proved simply  as  follows. 
From (\ref{3.2}) follows
$$
h_{\{r\}j_{1}...j_{n}u} = f_{rj_{1}...j_{n}u}
$$
\noindent and now (\ref{3.6}) is a consequence of (\ref{3.2}).
\par
Finally we notice that (\ref{3.4}) and (\ref{3.5}) are also consequences
of (\ref{3.6}) and (\ref{3.7}), i.e. \quad (\ref{3.4})\quad  $\Leftrightarrow$  \quad (\ref{3.6}) \quad   and \quad  (\ref{3.5})  
\quad $\Leftrightarrow$ \quad  (\ref{3.7}).\qquad \qquad $\Box$ 
 
Now we are ready to give the main theorem.

{\bf Theorem 3.2.} (i) {\it The  integrability  conditions  for  
the system} (\ref{3.1}) {\it for arbitrary initial conditions} 
$y_{i}(0, \cdots ,0)=C_{i}$, $1\le i\le n$, {\it are}
\begin{equation}
R_{\alpha \beta uv}  \equiv  0 \qquad {\it i}.{\it e}. \qquad 
R_{\{r\}\beta uv} \equiv  0.
\label{3.8}
\end{equation}
\par
\noindent (ii) {\it If the integrability conditions}  (\ref{3.8}) {\it are  
satisfied,
then there exist functions} 
$$P^{<w_{1},...,w_{k}>}_{i_{1}...i_{n}j_{1}...j_{n}} (x_{1},\cdots 
,x_{k}) 
,w_{1},\cdots 
,w_{n}\in \mathbb{ N}_{0}, i_{1},\cdots 
,i_{n},j_{1},\cdots 
,j_{n}\in \mathbb{ Z}$$
\par
\noindent {\it in a neighborhood of} $(0,\cdots ,0)$ {\it such that}
\begin{equation}
P^{<0,...,0>}_{i_{1}...i_{n}j_{1}...j_{n}} = \delta _{i_{1}j_{1}} \delta 
_{i_{2}j_{2}} ... \delta _{i_{n}j_{n}},
\label{3.9a}
\end{equation}
$$
P^{<w_{1},...,w_{u}+1,...,w_{k}>}_{i_{1}...i_{n}j_{1}...j_{n}} = {\partial 
\over \partial x_{u}} P^{<w_{1},...,w_{k}>}_{i_{1}...i_{n}j_{1}...j_{n}} +
$$
\begin{equation}
+\sum^{}_{t_{1},...,t_{n}\in \mathbb{ Z}} \Bigl( 
\sum^{n}_{s=1} i_{s}f
_{s(t_{1}-i_{1})...(t_{s}-i_{s}+1)...(t_{n}-i_{n})u} 
\Bigr) P^{<w_{1},...,w_{k}>}_{t_{1}...t_{n}j_{1}...j_{n}}. 
\label{3.9b}
\end{equation}
\noindent {\it If} $(C_{1}, \cdots ,C_{n})\in W$, 
{\it then the solution of} (\ref{3.1}) {\it in 
a neighborhood of} $(0, \cdots ,0)$ {\it is given by} 
$$
y_{1} = \sum^{}_{w_{1},...,w_{k}\in \mathbb{ N_{0}}} \Bigl[ 
{(-x_{1})^{w_{1}}\over w_{1}!}\cdots 
{(-x_{k})^{w_{k}}\over w_{k}!} \sum^{}_{j_{1},...,j_{n}\in \mathbb{ Z}} 
P^{<w_{1},...,w_{k}>}_{10...0j_{1}...j_{n}}C^{j_{1}}_{1} \cdot 
C^{j_{2}}_{2} \cdots 
C^{j_{n}}_{n}\Bigr] , 
$$
\begin{equation}
y_{2} = \sum^{}_{w_{1},...,w_{k}\in \mathbb{ N_{0}}} \Bigl[ 
{(-x_{1})^{w_{1}}\over w_{1}!}\cdots 
{(-x_{k})^{w_{k}}\over w_{k}!} 
\sum^{}_{j_{1},...,j_{n}\in \mathbb{ Z}} 
P^{<w_{1},...,w_{k}>}_{01...0j_{1}...j_{n}}C^{j_{1}}_{1}\cdot 
C^{j_{2}}_{2}\cdots 
C^{j_{n}}_{n} \Bigr] , 
\label{3.10} 
\end{equation}
$$
. . . . . . . . . . . . . . . . . . . . . . . .  
$$
$$
y_{n} = \sum^{}_{w_{1},...,w_{k}\in \mathbb{ N_{0}}} \Bigl[ 
{(-x_{1})^{w_{1}}\over w_{1}!}\cdots 
{(-x_{k})^{w_{k}}\over w_{k}!} 
\sum^{}_{j_{1},...,j_{n}\in \mathbb{ Z}} 
P^{<w_{1},...,w_{k}>}_{0...01j_{1}...j_{n}}C^{j_{1}}_{1}\cdot 
C^{j_{2}}_{2}\cdots 
C^{j_{n}}_{n} \Bigr] . 
$$
\noindent {\it This solution is unique with the given 
initial conditions in a neighborhood of} $(0, \cdots ,0)$.
\par
Proof. Let us introduce the following functions
$$
y_{\alpha } = y_{i_{1}i_{2}...i_{n}} = 
y^{i_{1}}_{1}\cdot y^{i_{2}}_{2} \cdots 
 y^{i_{n}}_{n}, \qquad (i_{1},\cdots 
,i_{n}\in \mathbb{ Z})
$$
\noindent such that $y_{1}=y_{\{1\}},\cdots 
,y_{n}=y_{\{n\}}$. These functions satisfy
$$
\partial y_{\{r\}}/\partial x_{u} + \sum^{}_{\alpha } 
f_{r\alpha u}y_{\alpha } = 0\qquad (1\le r\le n, 1\le u\le k)
$$
\noindent and hence
$$
\partial y_{\alpha }/\partial x_{u} = {\partial \over \partial 
x_{u}} (y^{i_{1}}_{1}\cdot y^{i_{2}}_{2} \cdots 
 y^{i_{n}}_{n}) 
$$
$$
= i_{1}y_{\alpha -\{1\}}\partial y_{1}/\partial x_{u} + \cdots 
 + i_{n}y_{\alpha -\{n\}}\partial y_{n}/\partial x_{u} 
$$
$$
= i_{1}y_{\alpha -\{1\}} \Bigl( 
-\sum^{}_{\beta }f _{1\beta u}y_{\beta } \Bigr) + \cdots + 
i_{n}y_{\alpha -\{n\}} \Bigl( 
-\sum^{}_{\beta } f_{n\beta u}y_{\beta } \Bigr) 
$$
$$
= -i_{1} \sum^{}_{\beta }
f_{1\beta u}y_{\alpha +\beta -\{1\}} - \cdots -
i_{n} \sum^{}_{\beta } 
f_{n\beta u}y_{\alpha +\beta -\{n\}}
$$
$$
= -\sum^{n}_{s=1}i_{s}\sum^{}_{\beta } 
f_{s\beta u}y_{\alpha + \beta -\{s\}}
$$
$$
= -\sum^{n}_{s=1}i_{s}\sum^{}_{\gamma } 
f_{s(\gamma -\alpha +\{s\})u}y_{\gamma },\qquad \hbox{ i.e.}
$$
\begin{equation}
{\partial \over \partial x_{u}} y_{\alpha } 
+\sum^{}_{\gamma }
h_{\alpha \gamma u}y_{\gamma } = 0
\label{3.11}
\end{equation}
\par
\noindent for $\alpha \in \mathbb{ Z}^{n}$, 
$u\in \{1,\cdots ,k\}$. 
\par
Thus we obtain that the  system  (\ref{3.1})  induces  the  system
(\ref{3.11}). Conversely also holds, i.e. one can prove that 
if the functions $f_{r\alpha u}$ are given, and
\par
(i) the system (\ref{3.11}) is satisfied, where $h_{\alpha \gamma 
u}$ are defined by (\ref{3.2}), 
\par
 (ii) $y_{i_{1}...i_{n}}(0,\cdots 
,0) = C^{i_{1}}_{1}\cdot C^{i_{2}}_{2}\cdots  C^{i_{n}}_{n}$ 
\qquad ($C_{i}$ are 
constants),
\par
\noindent then the system (\ref{3.1}) is satisfied, 
where $y_{r}=y_{\{r\}}$ for 
$1\le r\le n$.
\par
Similarly to the integrability conditions for the system 
(\ref{2.1}), the integrability conditions for the homogeneous linear 
system (\ref{3.11}) are given by $R_{\alpha 
\beta uv}\equiv 0$, i.e. 
$R_{\{s\}\beta uv}\equiv 0,$ because (\ref{3.7}) is
satisfied. Hence, the integrability  conditions  of  (\ref{3.1}) 
are given by (\ref{3.8}), 
and (i) is proved. 
\par 
Similarly to the proof of theorem 2.1, if 
the integrability conditions  (\ref{3.8})  are  satisfied,  then  there
exist functions
$$
P^{<w_{1},...,w_{k}>}_{\alpha \beta}(x_{1},\cdots 
,x_{k}),\qquad w_{1},\cdots 
,w_{k}\in \mathbb{ N}_{0}, \alpha , \beta \in \mathbb{ Z}^{n}
$$
\noindent such that (\ref{3.9a}) and (\ref{3.9b}) are satisfied. In order to 
prove that they are well defined, the convergence in (\ref{3.9b}) 
should be verified. It is easy to prove from (\ref{3.9a}) and 
(\ref{3.9b}) that for each $w_{1},\cdots ,w_{k}\in \mathbb{ N}_{0}$ 
the series 
$$\sum^{}_{i_{1},...,i_{n}\in \mathbb{ Z}} 
P_{i_{1}\cdots i_{n}j_{1}\cdots j_{n}}^{<w_{1},\cdots ,w_{k}>} 
z_{1}^{j_{1}-i_{1}}\cdot z_{2}^{j_{2}-i_{2}}\cdots 
z_{n}^{j_{n}-i_{n}}$$
uniformly converge for $(z_{1},\cdots ,z_{n})$ 
in a closed subset of W. The proof is by 
induction of $w_{1},\cdots ,w_{k}$ and it is analogous 
to the proof of the convergence of $\sum_{\gamma }
h_{\gamma \alpha v}h_{\beta \gamma u}$. The convergence in 
(\ref{3.9b}) follows simultaneously from here. 
Further by induction of $w_{1},\cdots ,w_{k}$ it also verifies the 
uniform convergence of 
$$\sum^{}_{j_{1},...,j_{n}\in \mathbb{ Z}} 
P_{i_{1}\cdots i_{n}j_{1}\cdots j_{n}}^{<w_{1},\cdots ,w_{k}>} 
z_{1}^{j_{1}}\cdot z_{2}^{j_{2}}\cdots 
z_{n}^{j_{n}}$$
\noindent for $(z_{1},\cdots ,z_{n})$ in a closed 
subset of W. Moreover, for fixed $i_{1},\cdots ,i_{n}$ there exist 
constants $M_{i_{1}\cdots i_{n}}$ in a neighborhood 
of the considered point, such that 
$$\mid \sum^{}_{j_{1},...,j_{n}\in \mathbb{ Z}} 
P_{i_{1}\cdots i_{n}j_{1}\cdots j_{n}}^{<w_{1},\cdots ,w_{k}>} 
C_{1}^{j_{1}}\cdots C_{n}^{j_{n}}\cdot 
{1\over w_{1}!\cdots w_{k}!} \mid \le M_{i_{1}\cdots i_{n}}
$$
\noindent for arbitrary $w_{1},\cdots ,w_{n}\in \mathbb{ N}_{0}$ 
and $(C_{1},\cdots ,C_{n})\in W$. The proof follows 
from a formula analogous to (\ref{2.9}). Indeed, using the same 
notations as in the proof of the theorem 2.1, by induction 
of $w_{1},\cdots ,w_{k}\in \mathbb{ N}_{0}$, it verifies that 
$$
\mid \sum ^{}_{\beta} Q_{\alpha \beta}^{<w_{1},\cdots ,w_{k}>} 
\sum ^{}_{\gamma } M_{\beta \gamma} 
C_{1}^{j_{1}}\cdots C_{n}^{j_{n}}\mid\le N_{\alpha} \qquad 
(\gamma =j_{1}\cdots j_{n})
$$
\noindent where $N_{\alpha }$ do not depend on $w_{1},\cdots ,
w_{k}$, and where 
$$
M_{\beta \gamma }=
\max _{a_{1},...,a_{k}\in \{ 0,1\} }
\mid {\partial ^{w_{1}-m_{1}+...+w_{k}-m_{k}}\over \partial 
x^{w_{1}-m_{1}}_{1}...\partial x^{w_{k}-m_{k}}_{k}} P^{<a_{1},\cdots 
,a_{k}>}_{\beta \gamma }\mid \cdot 
$$
$$
\cdot {\rho ^{w_{1}-m_{1}+...+w_{k}-m_{k}}\over 
(w_{1}-m_{1})!\cdots (w_{k}-m_{k})!}, 
$$
\noindent according to the Cauchy integral formula. 
\par 
Similarly to the proof of theorem 2.1, one can verify that 
the solution of (\ref{3.1}) with $y_{i_{1}\cdots i_{n}}(0,\cdots ,0)=
C^{i_{1}}_{1}\cdot C^{i_{2}}_{2} \cdots C^{i_{n}}_{n}$ 
is given by 
$$
y_{\alpha }=\sum^{}_{w_{1},...,w_{k}\in \mathbb{N}_{0}} \Bigl[ 
{(-x_{1})^{w_{1}}\over w_{1}!}\cdots 
{(-x_{k})^{w_{k}}\over w_{k}!}
\sum^{}_{\beta }
P^{<w_{1},...,w_{k}>}_{\alpha \beta} C^{j_{1}}_{1}C^{j_{2}}_{2}
\cdots 
C^{j_{n}}_{n} \Bigr]
$$
where $\beta =j_{1}\cdots j_{n}$. 
Its convergence follows 
from the previous discussion. 
Specially, if $\alpha \in \{1\},\cdots ,\alpha \in \{n\}$
we obtain the required solution (\ref{3.10}).      
\par
Each solution of (\ref{3.1})  is  analytical  function.  On  the 
other hand, by successive  differentiation  of  (\ref{3.1}),  we 
notice that all partial derivatives of $y_{s}$ can be 
calculated uniquely at $(0, \cdots ,0)$. Hence (\ref{3.1}) does 
not have more than one solution in a neighborhood of $(0, 
\cdots ,0)$, and the obtained solution is unique. 
\qquad $\Box$ 
\par
{\bf Remark 2.} If one (or more) equations in the system (\ref{3.1}) is (are) omitted, then 
we may consider the whole system with $nk$ equations, where 
the corresponding functions $F_{ru}$ can be considered as parameters. These parametric functions 
appear in the integrability conditions and also in the solution, and it leads to the parametric solution of the initial system.

\end{document}